\documentclass[11pt]{amsart}
\usepackage{amssymb,amsmath}
\usepackage{graphicx}
\usepackage[all]{xy}

\title[T-homotopy and refinement of observation (I)]{T-homotopy and refinement of observation (I) : Introduction}
\author[P. Gaucher]{Philippe Gaucher}
\address{Preuves Programmes et Syst{\`e}mes\\ Universit{\'e} Paris 7--Denis Diderot\\
Case 7014\\2 Place Jussieu\\ 75251 PARIS Cedex 05\\ France}
\email{gaucher@pps.jussieu.fr}
\urladdr{http://www.pps.jussieu.fr/{\~{}}gaucher/}

\newcommand{\C}{\mathcal{C}}

\newcommand{\de}{\partial}
\newcommand{\p}\times
\renewcommand{\vec}{\overrightarrow}
\renewcommand{\P}{\mathbb{P}}

\newcommand{\be}{\begin{equation}}
\newcommand{\ee}{\end{equation}}
\newcommand{\bea}{\begin{eqnarray}}
\newcommand{\eea}{\end{eqnarray}}
\newcommand{\beas}{\begin{eqnarray*}}
\newcommand{\eeas}{\end{eqnarray*}}

\newtheorem{thm}{Theorem}[section]
\newtheorem{prop}[thm]{Proposition}

\newtheorem{rem}[thm]{Remark}

\newtheorem{defn}[thm]{Definition}

\newcommand{\bd}{\begin{defn}}
\newcommand{\ed}{\end{defn}}
\newcommand{\bcd}{\begin{defn}}
\newcommand{\ecd}{\end{defn}}
\newcommand{\bex}{\begin{exmp}}
\newcommand{\eex}{\end{exmp}}
\newcommand{\bp}{\begin{prop}}
\newcommand{\ep}{\end{prop}}
\newcommand{\bth}{\begin{thm}}
\renewcommand{\eth}{\end{thm}}
\newcommand{\br}{\begin{rem}}
\newcommand{\er}{\end{rem}}
\newcommand{\bpf}{\begin{proof}}
\newcommand{\epf}{\end{proof}}

\newcommand{\fl}[1]{\ar@{->}[l]_{#1}}
\newcommand{\fr}[1]{\ar@{->}[r]^-{#1}}
\newcommand{\fd}[1]{\ar@{->}[d]_{#1}}
\newcommand{\fu}[1]{\ar@{->}[u]^{#1}}
\newcommand{\f}[2]{\ar@{->}[#1]|{#2}}
\newcommand{\ff}[2]{\ar@2{->}[#1]|{#2}}
\newcommand{\frr}[1]{\ar@{->}[rr]^{#1}}

\newcommand{\iso}{\cong}

\newcommand{\ot}{\otimes}

\newcommand{\vI}{\vec{I}}

\renewcommand{\leq}{\leqslant}
\renewcommand{\geq}{\geqslant}

\def\cartesien{%
  \ar@{-}[]+R+<6pt,-2pt>;[]+RD+<6pt,-6pt>%
  \ar@{-}[]+D+<2pt,-6pt>;[]+RD+<6pt,-6pt>%
}
\def\cocartesien{%
  \ar@{-}[]+L+<-6pt,+2pt>;[]+LU+<-6pt,+6pt>%
  \ar@{-}[]+U+<-2pt,+6pt>;[]+LU+<-6pt,+6pt>%
}

\newcommand{\brm}[1]{\rm{\mathbf{#1}}}

\renewcommand{\top}{{\brm{Top}}}

\newcommand{\dtop}{{\brm{Flow}}}

\newcommand{\glob}{{\rm{Glob}}}

\makeatletter

\def\varholim@#1#2{%
  \vtop{\m@th\ialign{##\cr
    \hfil$#1\operator@font holim$\hfil\cr
    \noalign{\nointerlineskip\kern1.5\ex@}#2\cr
    \noalign{\nointerlineskip\kern-\ex@}\cr}}%
}
\def\holimproj{%
  \mathop{\mathpalette\varholim@{\leftarrowfill@\textstyle}}\nmlimits@
}
\def\holiminj{%
  \mathop{\mathpalette\varholim@{\rightarrowfill@\textstyle}}\nmlimits@
}

\makeatother

\DeclareMathOperator{\cell}{{\brm{cell}}}
\DeclareMathOperator{\cof}{{\brm{cof}}}
\DeclareMathOperator{\inj}{{\brm{inj}}}
\newcommand{\hda}{{\cell(\dtop)}}

\begin{document}

\begin{abstract} 
This paper is the extended introduction of a series of papers about
modelling T-homotopy by refinement of observation.  The notion of
T-homotopy equivalence is discussed. A new one is proposed and its
behaviour with respect to other construction in dihomotopy theory is
explained.
\end{abstract}

\maketitle


\section{About deformations of HDA}

The main feature of the two algebraic topological models of
\textit{higher dimensional automata} (or HDA) 
introduced in \cite{diCW} and in \cite{model3} is to provide a
framework for modelling continuous deformations of HDA corresponding
to subdivision or refinement of observation. Globular complexes and
flows are specially designed to model the
\textit{weak S-homotopy equivalences} (the spatial  deformations) and 
the \textit{T-homotopy equivalences} (the temporal deformations). The
first descriptions of spatial deformation and of temporal deformation
dates back from the informal and conjectural paper
\cite{ConcuToAlgTopo}.

Let us now explain a little bit what the spatial and temporal
deformations consist of before presenting the results. The
computer-scientific and geometric explanations of
\cite{diCW} must of course be preferred for a deeper understanding.

\begin{figure}
\begin{center}
\includegraphics[width=6cm]{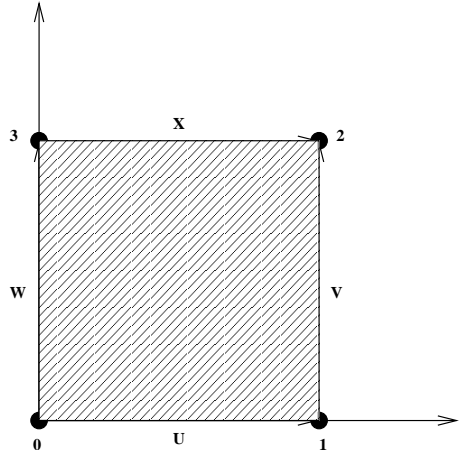}
\end{center}
\caption{Two concurrent processes}
\label{carreplein}
\end{figure}

In dihomotopy theory, processes running concurrently cannot be
distinguished by any observation. For instance in
Figure~\ref{carreplein}, each axis of coordinates represents one
process and the two processes are running concurrently. The
corresponding geometric shape is a full $2$-cube. This example
corresponds to the flow $\vec{C}_2$ defined as follows:
\begin{itemize}
\item Let us introduce the flow $\de\vec{C}_2$ defined 
by $(\de\vec{C}_2)^0=\{0,1,2,3\}$, $\P_{0,1}\de\vec{C}_2=\{U\}$,
$\P_{1,2}\de\vec{C}_2=\{V\}$, $\P_{0,3}\de\vec{C}_2=\{W\}$,
$\P_{3,2}\de\vec{C}_2=\{X\}$. The flow $\de\vec{C}_2$ corresponds 
to an empty square, where the execution paths $U*V$ and $W*X$ 
are \textit{not} running concurrently.
\item Then consider the pushout diagram 
\[
\xymatrix{
\glob(\mathbf{S}^{0})\fd{}\fr{q}& \de\vec{C}_2\fd{}\\
\glob(\mathbf{D}^1)\fr{} & \vec{C}_2 \cocartesien}
\]
with $q(\mathbf{S}^{0})=\{U*V,W*X\}$ (the globe functor $\glob(-)$ is
defined below). The presence of $\glob(\mathbf{D}^1)$ creates a
\textit{S-homotopy} between the execution paths $U*V$ and $W*X$,
modelling this way the concurrency.
\end{itemize}
It does not matter for $\P_{0,2}\vec{C}_2$ to be homeomorphic to
$\mathbf{D}^1$ or only homotopy equivalent to $\mathbf{D}^1$, or even
only weakly homotopy equivalent to $\mathbf{D}^1$. The only fact that
matters is that the topological space $\P_{0,2}\vec{C}_2$ be weakly
contractible. Indeed, a hole like in the flow $\de \vec{C}_2$ (the
space $\P_{0,2}\de\vec{C}_2$ is the discrete space $\{U*V,W*X\}$ )
means that the execution paths $U*V$ and $W*X$ are not running
concurrently, and therefore that they are distinguishable by
observation. This kind of identification is well taken into account by
the notion of weak S-homotopy equivalence.  This notion is introduced
in \cite{diCW} in the framework of globular complexes, in
\cite{model3} in the framework of flows and it is proved that these
two notions are equivalent in \cite{model2}.

\begin{figure}
\begin{center}
\includegraphics[width=6cm]{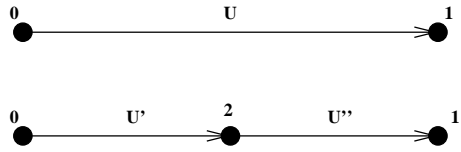}
\end{center}
\caption{The most simple example of T-homotopy equivalence}
\label{ex1}
\end{figure}

In dihomotopy theory, it is also required to obtain descriptions of
HDA which are invariant by refinement of observation. The simplest
example of refinement of observation is represented in
Figure~\ref{ex1}, in which the directed segment $U$ is divided in two
directed segments $U'$ and $U''$.  This kind of identification is well
taken into account by the notion of T-homotopy equivalence. This
notion is introduced in \cite{diCW} in the framework of globular
complexes, and in \cite{model2} in the framework of flows. The latter
paper also proves that the two notions are equivalent. In the case of
Figure~\ref{ex1}, the T-homotopy equivalence is the unique morphism of
flows sending $U$ to $U'*U''$.

Each weak S-homotopy equivalence as well as each T-homotopy
equivalence preserves as above the initial states and the final states
of a flow. More generally, any good notion of dihomotopy equivalence
must preserve the \textit{branching and merging homology theories}
introduced in \cite{exbranch}. This paradigm dates from the beginning
of dihomotopy theory: a dihomotopy equivalence must not change the
topological configuration of branching and merging areas of execution
paths \cite{survol}.  It is also clear that any good notion of
dihomotopy equivalence must preserve the \textit{underlying homotopy
type}, that is the topological space, defined only up to weak homotopy
equivalence, obtained after removing the time flow. In the case of
Figure~\ref{carreplein} and Figure~\ref{ex1}, this underlying homotopy
type is the one of the point.

\section{Prerequisites and notations}

The initial object (resp. the terminal object) of a category $\C$, if
it exists, is denoted by $\varnothing$ (resp. $\mathbf{1}$).

Let $\C$ be a cocomplete category.  If $I$ is a set of morphisms of
$\C$, then the class of morphisms of $\C$ that satisfy the RLP
(\textit{right lifting property}) with respect to any morphism of $I$
is denoted by $\inj(I)$ and the class of morphisms of $\C$ that are
transfinite compositions of pushouts of elements of $I$ is denoted by
$\cell(I)$. Denote by $\cof(I)$ the class of morphisms of $\C$ that
satisfy the LLP (\textit{left lifting property}) with respect to any
morphism of $\inj(I)$.  It is a purely categorical fact that
$\cell(I)\subset \cof(I)$. Moreover, every morphism of $\cof(K)$ is a
retract of a morphism of $\cell(K)$ as soon as the domains of $K$ are
small relative to $\cell(K)$ (\cite{MR99h:55031} Corollary~2.1.15). An
element of $\cell(I)$ is called a \textit{relative $I$-cell
complex}. If $X$ is an object of $\C$, and if the canonical morphism
$\varnothing\longrightarrow X$ is a relative $I$-cell complex, one
says that $X$ is a \textit{$I$-cell complex}.

Let $\C$ be a cocomplete category with a distinguished set of
morphisms $I$. Then let $\cell(\C,I)$ be the full subcategory of $\C$
consisting of the objects $X$ of $\C$ such that the canonical morphism
$\varnothing\longrightarrow X$ is an object of $\cell(I)$. In other
terms, $\cell(\C,I)=(\varnothing\!\downarrow \! \C) \cap \cell(I)$.

Possible references for \textit{model categories} are
\cite{MR99h:55031}, \cite{ref_model2} and \cite{MR1361887}.  The 
original reference is \cite{MR36:6480} but Quillen's axiomatization is
not used in this paper. The axiomatization from Hovey's book is
preferred.  If $\mathcal{M}$ is a cofibrantly generated model category
with set of generating cofibrations $I$, let $\cell(\mathcal{M}) :=
\cell(\mathcal{M},I)$. A cofibrantly generated model structure
$\mathcal{M}$ comes with a \textit{cofibrant replacement functor}
$Q:\mathcal{M} \longrightarrow \cell(\mathcal{M})$.

A \textit{partially ordered set} $(P,\leq)$ (or \textit{poset}) is a
set equipped with a reflexive antisymmetric and transitive binary
relation $\leq$. A poset is \textit{locally finite} if for any
$(x,y)\in P\p P$, the set $[x,y]=\{z\in P,x\leq z\leq y\}$ is
finite. A poset $(P,\leq)$ is \textit{bounded} if there exist
$\widehat{0}\in P$ and $\widehat{1}\in P$ such that $P\subset
[\widehat{0},\widehat{1}]$ and such that $\widehat{0} \neq
\widehat{1}$. Let $\widehat{0}=\min P$ (the bottom element) and
$\widehat{1}=\max P$ (the top element).

The category $\top$ of \textit{compactly generated topological spaces}
(i.e. of weak Hausdorff $k$-spaces) is complete, cocomplete and
cartesian closed (more details for this kind of topological spaces in
\cite{MR90k:54001,MR2000h:55002}, the appendix of \cite{Ref_wH} and
also the preliminaries of \cite{model3}). For the sequel, any
topological space will be supposed to be compactly generated. A
\textit{compact space} is always Hausdorff.

The time flow of a higher dimensional automaton is encoded in an
object called a \textit{flow} \cite{model3}. A flow $X$ consists of a
set $X^0$ called the \textit{$0$-skeleton} and whose elements
correspond to the \textit{states} (or \textit{constant execution
paths}) of the higher dimensional automaton. For each pair of states
$(\alpha,\beta)\in X^0\p X^0$, there is a topological space
$\P_{\alpha,\beta}X$ whose elements correspond to the
\textit{(nonconstant) execution paths} of the higher dimensional
automaton \textit{beginning at} $\alpha$ and
\textit{ending at} $\beta$. If $x\in \P_{\alpha,\beta}X$ , let
$\alpha=s(x)$ and $\beta=t(x)$. For each triple
$(\alpha,\beta,\gamma)\in X^0\p X^0\p X^0$, there exists a continuous
map $*:\P_{\alpha,\beta}X\p \P_{\beta,\gamma}X\longrightarrow
\P_{\alpha,\gamma}X$ called the \textit{composition law} which is supposed 
to be associative in an obvious sense. The topological space $\P
X=\bigsqcup_{(\alpha,\beta)\in X^0\p X^0}\P_{\alpha,\beta}X$ is called
the \textit{path space} of $X$. The category of flows is denoted by
$\dtop$. A point $\alpha$ of $X^0$ such that there are no non-constant
execution paths ending to $\alpha$ (resp. starting from $\alpha$) is
called an \textit{initial state} (resp. a \textit{final state}). A
morphism of flows $f$ from $X$ to $Y$ consists of a set map $f^0:X^0
\longrightarrow Y^0$ and a continuous map $\P f: \P X \longrightarrow
\P Y$ preserving the structure. A flow is therefore ``almost'' a small 
category enriched in $\top$.

The category $\dtop$ is equipped with the unique model structure
such that \cite{model3}: 
\begin{itemize}
\item The weak equivalences are the \textit{weak S-homotopy equivalences}, 
i.e. the morphisms of flows $f:X\longrightarrow Y$ such that
$f^0:X^0\longrightarrow Y^0$ is a bijection and such that $\P f:\P
X\longrightarrow \P Y$ is a weak homotopy equivalence. 
\item The fibrations are the morphisms of flows
$f:X\longrightarrow Y$ such that $\P f:\P X\longrightarrow \P Y$ is a
Serre fibration. 
\end{itemize}
This model structure is cofibrantly generated. The set of generating
cofibrations is the set $I^{gl}_+=I^{gl}\cup \{R,C\}$ with
\[I^{gl}=\{\glob(\mathbf{S}^{n-1})\subset \glob(\mathbf{D}^{n}), n\geq
0\}\] where $\mathbf{D}^{n}$ is the $n$-dimensional disk, where
$\mathbf{S}^{n-1}$ is the $(n-1)$-dimensional sphere, where $R$ and
$C$ are the set maps $R:\{0,1\}\longrightarrow \{0\}$ and
$C:\varnothing\longrightarrow \{0\}$ and where for any topological
space $Z$, the flow $\glob(Z)$ is the flow defined by
$\glob(Z)^0=\{\widehat{0},\widehat{1}\}$, $\P \glob(Z)=Z$,
$s=\widehat{0}$ and $t=\widehat{1}$, and a trivial composition
law. The set of generating trivial cofibrations is
\[J^{gl}=\{\glob(\mathbf{D}^{n}\p\{0\})\subset
\glob(\mathbf{D}^{n}\p [0,1]), n\geq 0\}.\]

\begin{figure}
\begin{center}
\includegraphics[width=7cm]{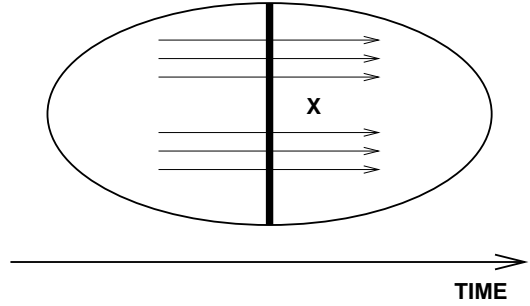}
\end{center}
\caption{Symbolic representation of
$\glob(X)$ for some topological space $X$} \label{exglob}
\end{figure}

\section{Why adding new T-homotopy equivalences ?}

\begin{figure}
\begin{center}
\includegraphics[width=6cm]{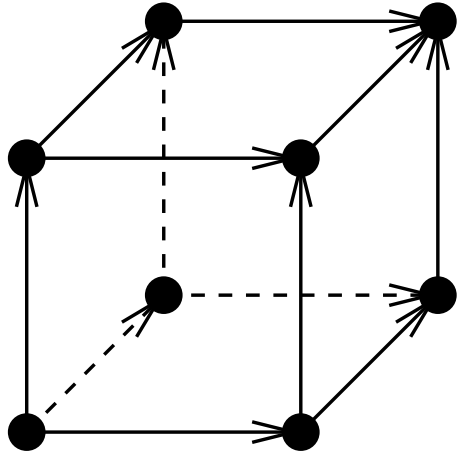}
\end{center}
\caption{The full $3$-cube}
\label{3cube}
\end{figure}

It turns out that the T-homotopy equivalences, as defined in
\cite{model2}, are the deformations which locally act like in
Figure~\ref{ex1}~\footnote{This fact was of course not known when
\cite{diCW} was being written down. The definition of T-homotopy
equivalence presented in that paper was based on the notion of
homeomorphism and it sounded so natural...}. So it becomes impossible
with this old definition to identify the directed segment of
Figure~\ref{ex1} with the full $3$-cube of Figure~\ref{3cube} by a
zig-zag sequence of weak S-homotopy and of T-homotopy equivalences
preserving the initial state and the final state of the $3$-cube since
any point of the $3$-cube is related to three distinct edges
(cf. Theorem~\ref{pasiso}). This contradicts the fact that concurrent
execution paths cannot be distinguished by observation. More
precisely, one has:

\bp\label{deftenseur}
Let $X$ and $Y$ be two flows. There exists a unique structure
of flows $X\ot Y$ on the set $X\p Y$ such that
\begin{enumerate}
\item $(X\ot Y)^0=X^0\p Y^0$
\item $\P (X\ot Y)= (\P X\p \P Y)\cup (X^0 \p \P Y) \cup (\P X\p Y^0)$
\item $s(x,y)=(s(x),s(y))$, $t(x,y)=(t(x),t(y))$, $(x,y)*(z,t)=(x*z,y*t)$.
\end{enumerate}
\ep

\bd The {\rm directed segment} $\vI$ is the flow $\glob(Z)$ with 
$Z=\{u\}$. \ed

\bd Let $n\geq 1$. The {\rm full $n$-cube} $\vec{C}_n$ is by definition 
the flow $Q(\vI^{\ot n})$, where $Q$ is the cofibrant replacement
functor. \ed

Notice that for $n\geq 2$, the flow $\vI^{\ot n}$ is not
cofibrant. Indeed, the composition law contains relations. For
instance, with $n=2$, one has $(\widehat{0},u)*(u,\widehat{1}) =
(u,\widehat{0}) * (\widehat{1},u)$

\bth \label{pasiso} 
Let $n\geq 3$. There does not exist any zig-zag sequence
\[
\xymatrix{
\vec{C}_n=X_0 \fr{f_0} & X_1 & X_2 \fl{f_1}\fr{f_2} & \dots &  X_{2n}=\vI \fl{f_{2n-1}}
}
\]
where each $X_i$ is an object of $\hda$ and where each morphism $f_i$
is either a S-homotopy equivalence~\footnote{Recall that a morphism
between two objects of $\hda$ is a weak S-homotopy equivalence if and
only if it is a S-homotopy equivalence.} or a T-homotopy
equivalence. \eth

In the statement of Theorem~\ref{pasiso}, we suppose that each flow
$X_i$ belongs to $\hda$ because T-homotopy is only defined between
this kind of flow.

\section{Full directed ball}

We need to generalize the notion of subdivision of the directed
segment $\vI$.

\bd 
A flow $X$ is {\rm loopless} if for every $\alpha\in X^0$, the space
$\P_{\alpha,\alpha}X$ is empty.
\ed

A flow $X$ is loopless if and only if the transitive closure of the
set $\{(\alpha,\beta)\in X^0\p X^0 \hbox{ such that }
\P_{\alpha,\beta}X \neq \varnothing\}$ induces a partial ordering on
$X^0$.

\bd 
A {\rm full directed ball} is a flow $\vec{D}$ such that:
\begin{itemize}
\item the $0$-skeleton $\vec{D}^0$ is finite
\item $\vec{D}$ has exactly one initial state $\widehat{0}$ and one final state $\widehat{1}$ with $\widehat{0} \neq \widehat{1}$
\item each state $\alpha$ of $\vec{D}^0$ is between $\widehat{0}$ and $\widehat{1}$, that is there exists an execution path from $\widehat{0}$ to $\alpha$, and another execution path from $\alpha$ to $\widehat{1}$
\item $\vec{D}$ is loopless
\item for any $(\alpha,\beta)\in \vec{D}^0\p \vec{D}^0$, the topological space 
$\P_{\alpha,\beta}\vec{D}$ is empty or weak\-ly contractible. 
\end{itemize}
\ed

Let $\vec{D}$ be a full directed ball. Then the set $\vec{D}^0$ can be
viewed as a finite bounded poset. Conversely, if $P$ is a finite
bounded poset, let us consider the \textit{flow} $F(P)$
\textit{associated to} $P$: it is of course defined as the unique flow 
(up to isomorphism) $F(P)$ such that $F(P)^0=P$ and
$\P_{\alpha,\beta}F(P)=\{u\}$ if $\alpha<\beta$ and
$\P_{\alpha,\beta}F(P)=\varnothing$ otherwise. Then $F(P)$ is a full
directed ball and for any full directed ball $\vec{D}$, the two flows
$\vec{D}$ and $F(\vec{D}^0)$ are weakly S-homotopy equivalent.

Let $\vec{E}$ be another full directed ball. Let
$f:\vec{D}\longrightarrow\vec{E}$ be a morphism of flows preserving
the initial and final states. Then $f$ induces a morphism of posets
from $\vec{D}^0$ to $\vec{E}^0$ such that $f(\min
\vec{D}^0)=\min \vec{E}^0$ and $f(\max \vec{D}^0)=\max \vec{E}^0$. Hence 
the following definition:

\bd 
Let $\mathcal{T}$ be the class of morphisms of posets
$f:P_1\longrightarrow P_2$ such that:
\begin{enumerate}
\item The posets $P_1$ and $P_2$ are finite and bounded. 
\item The morphism of posets $f:P_1 \longrightarrow P_2$ is one-to-one; 
in particular, if $x$ and $y$ are two elements of $P_1$ with $x<y$,
then $f(x)<f(y)$.
\item One has $f(\min P_1)=\min P_2$ and  $f(\max P_1)=\max P_2$.
\end{enumerate}
Then a {\rm generalized T-homotopy equivalence} is a morphism of
$\cof(\{Q(F(f)),f\in\mathcal{T}\})$ where $Q$ is the cofibrant
replacement functor of $\dtop$.
\ed

\begin{figure}
\begin{center}
\includegraphics[width=9cm]{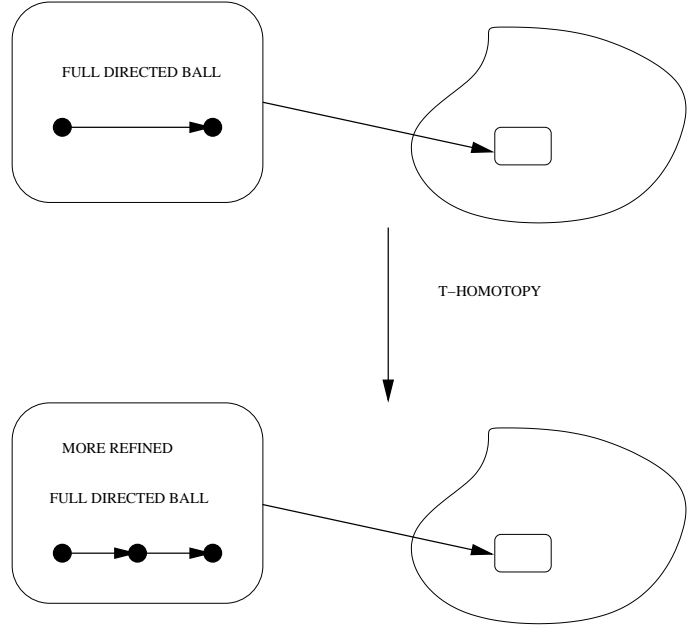}
\end{center}
\caption{Replacement of a full directed ball by a more refined one}
\label{ex2}
\end{figure}

A T-homotopy consists of locally replacing in a flow a full directed
ball by a more refined one (cf. Figure~\ref{ex2}).

\begin{figure}
\[
\xymatrix{
& A \fr{\leq} & B \ar@{->}[rd]^{\leq} & \\
\widehat{0} \ar@{->}[rd]^{\leq} \ar@{->}[ru]^{\leq} & & & \widehat{1}\\
&  C \ar@{->}[rru]^{\leq} & & }
\] 
\caption{Example of finite bounded poset}
\label{poset}
\end{figure}
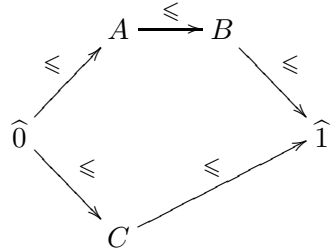

In a HDA, a $n$-transition, that is the concurrent execution of $n$
processes, is represented by the full $n$-cube $\vec{C}_n$. The
corresponding poset is the product poset
$\{\widehat{0}<\widehat{1}\}^n$. In particular, the poset corresponding 
to the full directed ball of Figure~\ref{3cube} is 
$\{\widehat{0}<\widehat{1}\}^3=\{\widehat{0}<\widehat{1}\}\p
\{\widehat{0}<\widehat{1}\}\p
\{\widehat{0}<\widehat{1}\}$.

The poset corresponding to Figure~\ref{carreplein} is the poset
$\{\widehat{0}<\widehat{1}\}^2=\{\widehat{0}<\widehat{1}\}\p
\{\widehat{0}<\widehat{1}\}$. If for instance $U$ is subdivided in two 
processes as in Figure~\ref{ex1}, the poset of the full directed ball
of Figure~\ref{carreplein} becomes equal to
$\{\widehat{0}<2<\widehat{1}\}\p
\{\widehat{0}<\widehat{1}\}$.

One has the isomorphism of flows $\vI^{\ot n} \iso
F(\{\widehat{0}<\widehat{1}\}^n)$ for every $n\geq 1$. The flow
$\vec{C}_n$ ($n\geq 1$) is identified to $\vI$ by the zig-zag sequence
of S-homotopy and generalized T-homotopy equivalences
\[
\xymatrix@1{
\vI & \fl{\simeq} Q(\vI) \ar@{->}[rr]^-{Q(F(g_n))} && Q(\vI^{\ot n}), }
\] 
where $g_n:\{\widehat{0} < \widehat{1}\} \longrightarrow \{\widehat{0} <
\widehat{1}\}^n \in \mathcal{T}$.

\section{Is this new definition well-behaved ?}

First of all, we must verify that each old T-homotopy equivalence as
defined in \cite{model2} will be a particular case of this new
definition. And indeed, one has:

\bth 
Let $X$ and $Y$ be two objects of $\hda$.  Let $f:X\longrightarrow Y$
be a T-homotopy equivalence as defined in \cite{model2}. Then $f$ can
be written as a composite $X\longrightarrow Z\longrightarrow Y$ where
$g:X\longrightarrow Z$ is a generalized T-homotopy equivalence and
where $h:Z\longrightarrow Y$ is a weak S-homotopy equivalence.
\eth

The two other tests consist of verifying that the branching and
merging homology theories \cite{exbranch}, as well as the underlying
homotopy type functor \cite{model2} are preserved with this new
definition of T-homotopy equivalence. And indeed, one has:

\bth 
Let $f:X\longrightarrow Y$ be a generalized T-homotopy equivalence.
Then for any $n\geq 0$, the morphisms of abelian groups
$H_n^-(f):H_n^-(X)\longrightarrow H_n^-(Y)$ and 
$H_n^+(f):H_n^+(X)\longrightarrow H_n^+(Y)$ are isomorphisms of groups
where $H_n^-$ (resp.$H_n^+$) is the $n$-th branching (resp. merging)
homology group. And the continuous map $|f|:|X|\longrightarrow |Y|$ is
a weak homotopy equivalence where $|X|$ denotes the underlying
homotopy type of the flow $X$.
\eth

\section{Conclusion}

This new definition of T-homotopy equivalence seems to be
well-behaved. It will hopefully have a longer lifetime than other ones
that the author proposed in the past. It is already known after
\cite{nonexistence} that it is impossible to construct a model
structure on $\dtop$ such that the weak equivalences are exactly the
weak S-homotopy equivalences and the generalized T-homotopy
equivalences. So new models of dihomotopy will be probably necessary
to understand the T-homotopy equivalences.

\end{document}